\documentclass[a4paper,reqno]{amsart}
\usepackage{amsmath,amssymb,eucal,bbold,amsthm}
\usepackage{array}
\usepackage[latin1]{inputenc}
\newcommand{\vac}{\vert 0\rangle}
\newcommand{\nvac}{a^\dagger(\ell_n)\ldots a^\dagger(\ell_1)\vac}
\newcommand{\tts}[1]{T_{1m_1}\ldots T_{1m_{#1}}}

\newcommand{\gS}{\mathfrak{S}}
\newcommand{\egy}{\mathcal{E}}
\newcommand{\integ}{\mathbb{Z}}
\newcommand{\comps}{\mathbb{C}}
\theoremstyle{plain}
\newtheorem{thm}{Theorem}
\newtheorem{prop}{Proposition}
\theoremstyle{definition}
\newtheorem{rmk}{Remark}

\begin{document}

\title{The energy operator for infinite statistics}
\author{Sonia Stanciu}
\address{Physikalisches Institut der Universität Bonn}
\thanks{BONN-HE-92-04} 
\date{February 1992}
\begin{abstract}
  We construct the energy operator for particles obeying infinite
  statistics defined by a $q$-deformation of the Heisenberg algebra.
\end{abstract}
\maketitle

The aim of this paper is to construct the energy operator for particles
obeying the so-called infinite statistics.  This topic was studied in
\cite{Zagier}, where a conjecture was formulated concerning the form of
the energy operator.  Our main result is a proof of this conjecture in
a slightly modified form (c.f.  Remark 1).

Infinite statistics is defined by the $q$-deformation of the
Heisenberg algebra, i.e., for all $k,\ell$,
\begin{align}
  \label{eq:COMM}
  & a(k) a^\dagger(\ell) -q a^\dagger(\ell) a(k) = \delta_{k,\ell}~, \\
  \label{eq:VAC}
  &a(k)\vac =0~,
\end{align}
where $\{a(k)\}$ represents a set of annihilation operators,
$a^\dagger(k)$ being their adjoints.  A state $x_{\bf k}$ is uniquely
specified by an ordered $n$-tuple $\boldsymbol{k} \equiv ( k_{1},
\ldots, k_{n})$, where we assume for simplicity $k_i\neq k_j,$ for
all $i,j$, so that
\begin{equation*}
  x_{\boldsymbol{k}} \equiv a^\dagger(k_n)\cdots a^\dagger(k_1)\vac~.  
\end{equation*}
Any other state which can be obtained by applying an arbitrary
monomial in the creation and annihilation operators on the vacuum can
be transformed using \eqref{eq:COMM} and \eqref{eq:VAC} into a linear
combination of the $x_{\boldsymbol{k}}$.

In order to show that these states build a Hilbert space
$\mathcal{H}(q)$, for a given $q$, one has to prove that the
corresponding inner product
$(x_{\boldsymbol{k}},x_{\boldsymbol{\ell}})$ yields a hermitian form,
i.e., that the infinite matrix
$A(q)=\{(x_{\boldsymbol{k}},x_{\boldsymbol{\ell}})\}$ is positive
definite.

Let us consider the inner product between two states
$x_{\boldsymbol{k}}=a^\dagger (k_m)\cdots a^\dagger (k_1)\vac$ and
$x_{\boldsymbol{\ell}}=a^\dagger (\ell_n)\cdots a^\dagger
(\ell_1)\vac$. This is easily seen to be zero unless $m=n$ and
$\boldsymbol{\ell}$ is a permutation of $\boldsymbol{k}$, so that
\begin{equation*}
  (x_{\pi (\boldsymbol{k})},x_{\boldsymbol{k}})= \langle0\vert
  a(k_{\pi (1)})\cdots a(k_{\pi (n)}) a^\dagger (k_{n})\cdots
  a^\dagger (k_{1})\vac =q^{I(\pi )}~;
\end{equation*}
or, consequently
\begin{equation*}
  (x_{\pi (\boldsymbol{k})},x_{\sigma
    (\boldsymbol{k})})=q^{I(\sigma^{-1}\pi)}~,
\end{equation*}
where $I(\pi )$ denotes the number of inversions of $\pi $, i.e.,
\begin{equation*}
  I(\pi )=\sum_{j=1}^n \text{Card}\{1\le i\le n \mid i<j
  \quad\text{and}\quad \pi(i)>\pi(j)\}~.
\end{equation*}
So the Hilbert space $\mathcal{H}(q)$  decomposes as an infinite
direct sum
\begin{equation*}
  \mathcal{H}(q)=\bigoplus_{n\geq 0}\bigoplus_{k}
  \mathcal{H}_{n,k}(q)~,
\end{equation*}
relative to which $A(q)$ decomposes as
\begin{equation*}
  A(q)=\oplus_{n\geq 0}\oplus_{k} A_{n,k}(q)~.
\end{equation*}
where $k$ denotes the \emph{unordered} $n$-tuple $\{k_{1},\ldots
,k_{n}\}$.  Each of the terms in these sums is finite dimensional and
$A_{n,k}$ does not actually depend on $k$
\begin{equation*}
  A_{n,k}(q)(\pi ,\sigma )\equiv A_{n}(q)(\pi ,\sigma
  )=q^{I(\sigma^{-1}\pi)}~,\qquad \forall \pi ,\sigma \in\gS_n~.
\end{equation*}
The positive definiteness of $A_n(q)$, for $-1<q<1$ has been proven by
D.~Zagier \cite{Zagier} who showed that
\begin{equation*}
  \det\,A_n(q)=\prod_{k=1}^n(1-q^{k^2+k})^{n!(n-k)\over k^2+k}~,
\end{equation*}
from which it follows that $A_n(q)$ is non-singular for all complex
numbers $q$ except the $N$-th roots of unity, for
$N=k^2+k=2,6,12,\ldots ,n^2+n$.  He also gave an explicit description
of $A_n^{-1}(q)$ and conjectured that
\begin{equation*}
  A_n^{-1}(q)\in \frac1{\Delta_n} M_{n!}(\integ[q])\quad,\quad
  \Delta_n \equiv (1-q^2)(1-q^6)\cdots(1-q^{n^2+n})~.
\end{equation*}

We begin by introducing the following notation.  If $\gS_n$ is the
group of permutations of $n$ elements we will denote by $T_{1k}$ the
particular elements which send $[1,2,\ldots ,n]$ to
$[k,1,\ldots,k-1,k+1,\ldots,n]$, i.e.,
\begin{equation*}
  T_{1k}(i)=
  \begin{cases}
    k,& \text{if $i=1$}\\
    i-1,& \text{if $1<i\le k$}\\
    i,& \text{if $k<i\le n$}
  \end{cases}
\end{equation*}
and by $\gS_{n,p}$ the following subsets of $\gS_n$:
\begin{equation*}
\gS_{n,p}=\{\sigma\in\gS_n,~{\rm with}\quad\sigma=
T_{1k_1}T_{1k_2}\cdots T_{1k_p},\quad 1< k_1<\cdots <k_p\le n\}~.
\end{equation*}

We will also consider particular elements in the group algebra 
\begin{equation*}
  \comps[\gS_n] = \left\{\sum_{\pi\in\gS_n} t_\pi \pi
    \quad\Bigg\vert\quad t_\pi \in \comps \right\}
\end{equation*}
of the form
\begin{equation}
  \label{eq:ALPHA}
  \alpha_n=\sum_{\rho\in \gS_n}A_n(\rho,1)\rho= \sum_{\rho\in
    \gS_n}q^{I(\rho)}\rho~,
\end{equation}
(We have considered $q$ fixed and omitted it from the notation.) Then, 
according to \cite{Zagier}, $\alpha_n$ is invertible in the group algebra 
if $\Delta_n\not= 0$ and $\alpha_n^{-1}$ will be given by
\begin{equation}
  \label{eq:ALPHAINV}
  \alpha_n^{-1}=\sum_{\rho\in
    \gS_n}A_n^{-1}(\rho,1)\rho~.
\end{equation}

Let $\egy$ be the energy operator of particles obeying infinite
statistics, defined by the commutation relation (1) in \cite{Zagier}.
$\egy$ acts on $\mathcal{H}(q)$ and each $x_{\boldsymbol{\ell}}$ is an
eigenvector of $\egy$ satisfying the eigenvalue equation
\begin{equation}
  \label{eq:EGY}
  \egy\nvac=\sum_{i=1}^n E(\ell_i)\nvac~,
\end{equation}
where $E(\ell_i)$ is the energy of a particle with momentum $\ell_i$.

\begin{thm}
  The energy operator $\egy$ has the form
  \begin{equation*}
    \egy =\sum_{n\ge 1}\egy_n~,
  \end{equation*}
  with
  \begin{equation}
    \label{eq:EGYN}
    \egy_n =\sum_{k_1,\ldots ,k_n}\sum_{\pi\in\gS_n}\sum_{i=1}^n
    c_i(q,\pi)E(k_{\pi(i)})a^\dagger(k_{\pi(n)})\cdots
    a^\dagger(k_{\pi(1)}) a(k_1)\cdots a(k_n),
  \end{equation}
  where the coefficients $c_i(q,\pi)$ are given by
  \begin{small}
    \begin{equation*}
      \sum_{\pi\in\gS_n}\sum_{i=1}^n c_i(q,\pi)X^{i-1}\pi=
      \alpha_n^{-1}\bigl(1-qXT_{12}\bigr)\bigl(1-q^2
      XT_{13}\bigr)\cdots \bigl(1-q^{n-1} XT_{1n}\bigr) \in
      \comps[X][\gS_n]
    \end{equation*}
  \end{small}
  or, explicitly,
  \begin{equation*}
    c_i(q,\pi)=(-1)^{i-1}\sum_{\tau\in\gS_{n,i-1}}
    A_n^{-1}(q)(\pi,\tau)A_n(q)(\tau,1),
  \end{equation*}
  for all $\pi\in\gS_n$ , $1\le i\le n$.
\end{thm}

\begin{rmk}
  The theorem agrees with Zagier's conjecture in \cite{Zagier} except
  that he has $E(k_i)$ instead of $E(k_{\pi(i)})$.  Thus the formulas
  agree if (and only if)
  \begin{equation*}
    c_i(q,\pi)=c_{\pi(i)}(q,\pi)~,
  \end{equation*}
  in for all $1\le i\le n$, $\pi\in\gS_n$.  This is true for $n\le 4$,
  but we do not know if it holds in general.
\end{rmk}

\begin{rmk}
  Although $\egy$ contains an infinite sum, when applied on a given
  $n$-particle state, only the first $n$ terms will give a nonzero
  contribution.
\end{rmk}

\begin{rmk}
  For $q=0$ this agrees with the results of O.~Greenberg \cite{Wally},
  who gave an expression for the energy operator of the form
  \begin{equation*}
    E=\sum_i \egy(i)n(i)~,
  \end{equation*}
  where the number operator $n(i)$ is given by
  \begin{equation*}
    n(i)=\sum_{s\ge 0}\sum_{k_1,\ldots, k_s}a^\dagger (k_1)\cdots a^\dagger(
    k_s) a^\dagger(i) a(i) a(k_s)\cdots a(k_1)~,
  \end{equation*}
  with obvious notation.  
\end{rmk}

To prove this theorem we need some preparation. We know from
\cite{Zagier} that the Hilbert space of states $\mathcal{H}(q)$ splits
into an infinite direct sum of finite dimensional blocks.  Each block
is determined by the \emph{unordered} $n$-tuple $\{k_1,\ldots,k_n\}$,
whereas a particular state in it is specified by an ordered version of
that $n$-tuple. In other words, we identify the Fock space states with
ordered sets $K=[k_1, \ldots,k_n]$.  For such a finite ordered set $A$
we denote by $s(A)$ and $l(A)$ respectively the smallest and the
largest element of $A$.  Ordered sets can be concatenated, e.g., if we
consider two disjoint ordered sets $A_1$ and $A_2$ we can form a new
ordered set $A_1\sqcup A_2$, such that if $a_i\in A_i$ then $a_1<a_2$.
Also, if $B$ is a subset of an ordered set $A$, one can form the
ordered set $A-B$. Moreover, we can invert the order of a given set,
the new one being denoted by $\overline{A}$.

The permutation group $\gS_n$ acts naturally on the ordered sets of
$n$ elements; and this action extends to an action of the group
algebra $\comps[\gS_n]$ on the vector space $\mathcal{L}$ of formal
linear combinations of such sets.  If $A$ is a given ordered set and
$\sigma\in\gS_n$ we define $I_A(\sigma A) = I(\sigma)$.
 
We conclude these general considerations by introducing a linear
evaluation map $\xi$ acting on $\mathcal{L}[X]$ and defined by
\begin{equation*}
  \xi\bigl(\sigma AX^{i-1}\bigr)=E\bigl((\sigma A)(i)\bigr)\sigma A~.
\end{equation*}

In order to be able to determine the coefficients $c_i (q,\pi)$ in
\eqref{eq:EGYN}, we have to understand how the energy operator $\egy$
and, in particular, each $\egy_p$ acts on an arbitrary state. For that
we will need two steps.

\begin{prop}
  The action of the $p$-particle term of the
  energy operator on given $n$-particle state $K$ is given by
  \begin{equation}
    \label{eq:EGYP}
    \egy_p K=\xi\Bigl(X^{n-p}\sum_{J\subset K\atop |J|=p} q^{I_K\bigl(
      (K-J)\sqcup J\bigr)}\bigl(K-J\bigr)\sqcup R_p(q,X) J\Bigr)~,
  \end{equation}
  where
  \begin{equation}
    R_p(q,X)=\alpha_p\sum_{\pi\in\gS_p}\sum_{i=1}^p
    c_i(q,\pi)X^{i-1}\pi~,
    \label{eq:REP}
  \end{equation}
  for all $1\le p\le n$.
\end{prop}

\begin{proof}
  To begin with, let us consider the case $p=n$.  We have
  (\cite{Zagier},\S 2)
  \begin{equation*}
    a(k_1)\cdots a(k_n)\nvac=\sum_{\sigma\in\gS_n}q^{I(\sigma)}
    \delta_{k_1 l_{\sigma(1)}}\cdots\delta_{k_n l_{\sigma(n)}}\vac~.
  \end{equation*}
  Thus, applying $\egy_n$ on an $n$-particle state we obtain
  \begin{equation*}
    \begin{aligned}
      \egy_n K &=\sum_{\sigma,\pi\in\gS_n}\sum_{i=1}^n q^{I(\sigma)}
      c_i(q,\pi)E((\sigma\pi K)(i))\sigma\pi K\\
      &=\xi\Bigl(\sum_{\sigma,\pi\in\gS_n}\sum_{i=1}^n q^{I(\sigma)}
      c_i(q,\pi)X^{i-1}\sigma\pi K\Bigr)\\
      &=\xi\Bigl(R_n(q,X)K\Bigr)~.
    \end{aligned}
  \end{equation*}
  
  We must now determine how a generic term $\egy_p$ acts on the
  $n$-particle state. Its action can be described in the following
  way: it chooses a subset $J\subset K$, $|J|=p$, such that the $p$
  annihilation operators of $\egy_{p}$ will ``contract'' with the $p$
  creation operators of $J$, leaving the remaining creation operators
  of $K$ in unaltered order {\it i.e.\/}, characterized by the set
  $(K-J)$. This yields a new $n$-particle state, characterized by the
  permutation $(K-J)\sqcup J$ multiplied by the numerical coefficient
  incurred in by repeated application of the commutation relation (1)
  in \cite{Zagier} and which is given by $q^{I_K\bigl( (K-J) \sqcup
    J\bigr)}$.  Clearly, $R_p(q,X)$ acts now on $J$ and, because the
  evaluation map $\xi$ is defined on the whole $n$-particle state, we
  have to shift the polynomial in $X$ by a common factor
  $X^{|K-J|}=X^{n-p}$ in order to obtain the correct energies.  Hence,
  it follows that
  \begin{equation*}
    \egy_p K=\xi\Bigl(X^{n-p}\sum_{J\subset K\atop |J|=p}
    q^{I_K\bigl( (K-J)\sqcup J\bigr)}\bigl(K-J\bigr)\sqcup
    R_p(q,X)J\Bigr)~.
  \end{equation*}
\end{proof}

\begin{prop}
  The action of the group ring element $R_p(q,X)$ on the ordered set
  $J$ is given by
  \begin{equation}
    \label{eq:WRONX}
    R_p(q,X)J=\sum_{L\subset J\atop s(J)\not\in
      L}q^{I_J\bigl(\overline{L} \sqcup(J-L)\bigr)}
    \bigl(\overline{L}\sqcup (J-L)\bigr)\bigl(-X\bigr)^{|L|}~.
  \end{equation}
\end{prop}

\begin{proof}
  We shall essentially show that \eqref{eq:WRONX} yields the correct
  energy operator,{\it i.e.}, that it satisfies the eigenvalue
  equation.  Therefore, we insert $R_p(q,X)$ in the expression for
  $\egy_p$ and we compute
  \begin{small}
    \begin{equation*}
      \egy_p K=\xi\Bigl(X^{n-p}\sum_{J\subset K\atop |J|=p}
      \sum_{L\subset J\atop s(J)\not\in L}
      q^{I_K\bigl((K-J)\sqcup J\bigr)+I_J\bigl(\overline{L}
        \sqcup(J-L)\bigr)} \bigl((K-J)\sqcup \overline{L}\sqcup
      (J-L)\bigr)\bigl(-X\bigr)^{|L|}\Bigr)~.
    \end{equation*}
  \end{small}
  But, obviously,
  \begin{equation*}
    I_K\bigl((K-J)\sqcup J\bigr)+I_J\bigl(\overline{L}
    \sqcup(J-L)\bigr)= I_K\bigl((K-J)\sqcup \overline{L}
    \sqcup(J-L)\bigr)~,
  \end{equation*}
  such that we obtain
  \begin{equation*}
    \egy_p K=\xi \Bigl(\sum_{J\subset K\atop |J|=p}
    \sum_{L\subset J\atop s(J)\not\in L}(-1)^{|L|}
    q^{I_K\bigl((K-J)\sqcup \overline{L} \sqcup(J-L)\bigr)}
    \bigl((K-J)\sqcup \overline{L}\sqcup (J-L)\bigr)X^{n-p+|L|}\Bigr)~.
  \end{equation*}
  
  For given $J$ and $L$, we consider those terms in the sum which are 
  characterized by $l(K-J)>l(L)$. Then the corresponding set can be 
  viewed in another way, namely,
  \begin{equation*}
    (K-J)\sqcup \overline{L}\sqcup (J-L)=\Bigl((K-J)-\{l(K-J)\}\Bigr)
    \sqcup\Bigl(\{l(K-J)\}\sqcup \overline{L}\Bigr)\sqcup (J-L)~,
  \end{equation*}
  having now $l\Bigl((K-J)-\{l(K-J)\}\Bigr)<l\Bigl(\{l(K-J)\}
  \sqcup \overline{L}\Bigr)$ and thus corresponding to another set which 
  contributes as well to the sum. As one can easily see, these two 
  terms will occur with identical coefficients but with 
  opposite signs and will therefore cancel.
  
  Thus, it only remains to discuss the case $L=\varnothing$. If
  $l(K-J)>s(J)$, then we can proceed analogously, writing
  \begin{equation*}
    (K-J)\sqcup \varnothing \sqcup J=\Bigl((K-J)-\{l(K-J)\}\Bigr)
    \sqcup\{l(K-J)\}\sqcup J~,
  \end{equation*}  
  such that we obtain the usual cancellation. But if $l(K-J)>s(J)$, then 
  $(K-J)\sqcup J= K$, and we finally obtain
  \begin{equation*}
    \begin{aligned}
      \egy K&=\sum_{p=1}^n \egy_p K\\
      &=\xi\Bigl(\sum_{p=1}^n\sum_{J\subset K\atop |J|=p}
      q^{I_K\bigl((K-J)\sqcup J\bigr)}(K-J)\sqcup J~X^{n-p}\Bigr)\\
      &=\sum_{i=1}^n E(K_i)K~.
    \end{aligned}
  \end{equation*}
\end{proof}

Now we are ready to prove the theorem stated at the very beginning.

\begin{proof}[Proof of the Theorem]
  Let us return now to the usual permutation language. One can 
  easily see that the permutations of the form $\overline{L}\sqcup (J-L)$, 
  with $|L|=s$ can be written as $\tts{s}$, with $1<m_1<\cdots <m_s\le n$, 
  so that 
  \begin{equation*}
    R_n(q,X)=\sum_{s=0}^{n-1}\sum_{1<m_1<\cdots <m_s\le n}(-1)^s 
    q^{(m_1 -1)+(m_2 -1)+\cdots+(m_s -1)}X^s\tts{s}~,
  \end{equation*}
  where we used the fact that $I(T_{1k})=k-1$.
  
  Now we only have to identify this expression obtained for $R_n(q,X)$ 
  with its definition \eqref{eq:REP}, and we obtain the desired result; that
  is, the generating function for the coefficients $c_i(q,\pi)$ is given
  by
  \begin{equation*}
    \sum_{\pi\in\gS_n}\sum_{i=1}^n c_i(q,\pi)X^{i-1}\pi=
    \alpha_n^{-1}\bigl(1-qXT_{12}\bigr)\bigl(1-q^2 XT_{13}\bigr)\cdots
    \bigl(1-q^{n-1} XT_{1n}\bigr)~,
  \end{equation*}
  with $\alpha_n$ given by \eqref{eq:ALPHA}.
  
  The coefficients $c_i(q,\pi)$ can be also given in another equivalent form. 
  Using \eqref{eq:ALPHAINV}, the right-hand side of the equation above can
  be written as 
  \begin{multline*}
        \Bigl(\sum_{\rho\in\gS_n}A_n^{-1}(\rho,1)\rho\Bigr)
        \Bigl(\sum_{i=1}^n\sum_{\pi\in\gS_{n,i-1}}(-1)^{i-1}
        A_n(\pi,1)X^{i-1}\pi\Bigr)\\
        =\sum_{i=1}^n (-1)^{i-1}X^{i-1} \sum_{{\sigma\in\gS_n}\atop
          {\pi\in\gS_{n,i-1}}} A_n^{-1}(\sigma,\pi)A_n(\pi,1)\sigma~,
  \end{multline*}
  where we made the substitution $\sigma=\rho\pi\in\gS_n$ and we used 
  the fact that $A^{-1}_n(\sigma,\pi)=A^{-1}_n(\sigma\pi^{-1},1)$.   
  Thus, identifying with the left-hand side, we obtain
  \begin{equation*}
    c_i(q,\sigma)=(-1)^{i-1}\sum_{\pi\in\gS_{n,i-1}}
    A_n^{-1}(\sigma,\pi)A_n(\pi,1)~.
  \end{equation*}
  
  It only remains to show that the solution obtained is unique. 
  First of all it is obvious that the form \eqref{eq:EGYN} of the 
  energy operator is the most general which can be assumed for 
  such a system, so that we only need to consider the possibility 
  of having another set of coefficients $c^*_i(q,\pi)$, such that the 
  corresponding $\egy^*$ yields the same eigenvalue equation. 
  Hence we must have $(\egy-\egy^*)K=0$, $\forall K$. If we 
  consider a $1$-particle state, we imediately obtain $\Delta 
  c_1(q,1)\equiv c_1(q,1)-c_1^*(q,1)=0$, for $n=1$. Assume now 
  $\Delta c_i(q,\pi)=0$, for all $1\le i\le p~,~\pi\in\gS_p$ 
  in all orders $1\le p\le n-1$. Then for an $n$-particle state 
  we will have $(\egy-\egy^*)K=(\egy_n-\egy_n^*)K=0$. But, on the 
  other hand
  \begin{equation*}
    (\egy_n-\egy_n^*)K=\Bigl(\sum_{\rho,\pi\in\gS_n}\sum_{i=1}^n 
    q^{I(\rho\pi^{-1})}\Delta c_i(q,\pi)E(\rho(i))\rho\Bigr)K~.
  \end{equation*}
  Taking into account the fact that $\rho$ and $E(\rho(i))$ are 
  linearly independent we get
  \begin{equation*}
    \sum_{\pi\in\gS_n}A_n(q)(\rho,\pi)\Delta c_i(q,\pi)=0\qquad 
    \forall\rho\in\gS_n~,1\le i\le n~.
  \end{equation*}
  Since $A_n(q)$ is invertible, it follows that $\Delta c_i(q,\pi)=0$. 
  Hence the energy operator is uniquely determined.
\end{proof}

\section*{Acknowledgments}

I would like to thank W.~Nahm for constant help and encouragement.  I
am grateful to D.~Zagier for suggesting me this problem and for a
careful reading of the manuscript. It is a pleasure to thank
V.~Rittenberg for bringing reference \cite{Wally} to my attention.  I
would also like to thank J.~M.~Figueroa-O'Farrill for the many useful
conversations and M.~Terhoeven and Th.~Wittlich for {\TeX}nical help.

\begingroup\raggedright%
\endgroup

\end{document}